\documentclass{amsart}

\newtheorem{theorem}{Theorem}[section]

\newtheorem{proposition}[theorem]{Proposition}
\newtheorem{corollary}[theorem]{Corollary}
\theoremstyle{definition}

\theoremstyle{remark}

\numberwithin{equation}{section}

\usepackage{tikz}
\usepackage{lmodern}
\usetikzlibrary{shapes,arrows}
\usepackage[psamsfonts]{amssymb}
\usepackage{amsmath}

\begin{document}

\title[Quintic Energy Critical Wave in 3D Cylindrical Convex Domains]{Well-Posedness for Quintic Energy Critical Wave in 3D Cylindrical Convex Domains}
\author{MEAS Len}
\address{Department of Mathematics, Royal University of Phnom Penh, Phnom Penh, Cambodia}
\email{meas.len@rupp.edu.kh}
    \thanks{This work was partially supported by a grant from the Niels Hendrik Abel Board.}


\subjclass[2020]{Primary 35L71, 35L05; Secondary 58J45, 42B25}


\keywords{Energy critical waves, cylindrical domains, dispersive estimates, Strichartz estimates}

\begin{abstract}
In this paper, we establish the well-posedness in energy space for the quintic energy critical wave inside a cylindrical convex domain $\Omega\subset\mathbb{R}^3$ with smooth boundary $\partial\Omega\neq\emptyset$. The key tools  to prove local well-posedness are the dispersive estimates obtained in \cite{L,L1,L3} and the Strichartz estimates in \cite{L2}. We point out that our result on the local and global existence of the solution to the wave equation in the cylindrical domain setting interpolates between that of in Euclidean space $\mathbb{R}^3$ (see \cite{MG}) and in any bounded domains in $\mathbb{R}^3$ (see \cite{BLP}). Moreover, the result of the Strichartz estimates in our setting is strong enough when combined with the arguments in \cite{BLP,SS2} so that  we can extend local to global well-posedness.
\end{abstract}

\maketitle



\section{Introduction}
Let $\Omega=\{x\geq0,(y,z)\in\mathbb{R}^2\}\subset\mathbb{R}^3$ be a convex domain with smooth boundary $\partial\Omega=\{x=0\}$ 
 and $\Delta=\partial_x^2+(1+x)\partial_y^2+\partial_z^2$ the Laplacian acting on functions with Dirichlet boundary condition. We consider the following energy critical semilinear wave equation
\begin{equation}\label{NLW}
\begin{aligned}
&(\partial_t^2-\Delta)u+u^5=0\quad \text{in } \mathbb{R}_t\times\Omega,\\
&u_{|t=0}=u_0,\,\,\, \partial_t u_{|t=0}=u_1,\,\,\, u_{|x=0}=0,
\end{aligned}
\end{equation}
 In this work, we focus on the questions of well-posedness of \eqref{NLW} for the initial condition \[(u,\partial_t u)_{\vert t=0}=(u_0,u_1)\in H_0^1(\Omega)\times L^2(\Omega).\]
 
 The Riemannian manifold $(\Omega,\Delta)$ with Laplacian $\Delta=\partial_x^2+(1+x)\partial_y^2+\partial_z^2$ can be locally seen as a cylindrical domain in $\mathbb{R}^3$ by taking cylindrical coordinates $(r,\theta,z)$, where we set
$r=1-\frac{x}{2},\theta=y,$ and $z=z$ (see Remark 1.1\cite{ILP}).  In our case of cylindrical domain, the nonnegative radius of curvature depends on the incident angle and vanishes in some directions, and the boundary is convex with zero curvature along the axis of the cylinder. We point out here that our domain interpolates between the Euclidean space $\mathbb{R}^3$ and the bounded domain in $\mathbb{R}^3$. 

The motivation to consider this Laplacian $\Delta$ in our context comes from the Friedlander's model  domain of the half space $\Omega_D=\{(x,y)\mid x>0, y\in\mathbb{R}\}$ with Laplace operator given by
\[
\Delta_D=\partial_x^2+(1+x)\partial_y^2.
\]
We point out that when there is no $z$ variable in our Laplacian, the problem is reduced to the Friedlander's model (see \cite{ILP}). Moreover, the Laplacian $\Delta$ in our setting has a nice feature that allows explicit computations. Finally, let us identify a connection between the Laplacian $\Delta$ and the Laplace-Beltrami operator $\Delta_g$. For a metric $g=dx^2+(1+x)^{-1}dy^2+dz^2$, we see that the Laplace-Beltrami is given by $\Delta_g=(1+x)^{\frac{1}{2}}\partial_x(1+x)^{-\frac{1}{2}}\partial_x+(1+x)\partial_y^2+\partial_z^2,$ which is a self-adjoint operator with the volume form $\sqrt{\det g}dxdydz=(1+x)^{-\frac{1}{2}}dxdydz$. This implies that the difference $\Delta_g-\Delta=-\frac{1}{2}(1+x)^{-1}\partial_x$ is the first order differential operator. Our model uses instead the Laplace operator associated with the Dirichlet form $$\int |\nabla u|^2 dxdydz=\int(|\partial_x u|^2+(1+x)|\partial_y u|^2+|\partial_zu|^2)dxdydz.$$

We now turn to see some key properties of the solution to \eqref{NLW}. The solution $u$ of \eqref{NLW} possesses the invariant property under the dilation symmetry
 $$u(t,x)\mapsto u_\lambda(t,x):=\lambda^{\frac{1}{2}}u(\lambda t,\lambda x)$$ 
and
$$\partial_tu(t,x)\mapsto \partial_tu_\lambda(t,x):=\lambda^{\frac{3}{2}}\partial_tu(\lambda t,\lambda x).$$ 
 Recall that the $H_0^1(\Omega)\times L^2(\Omega)$ is invariant under this dilation symmetry. Indeed, one has
\[
\|(u_\lambda(0),\partial_tu_\lambda(0))\|_{H_0^1(\Omega)\times L^2(\Omega)}=\|(u(0),\partial_tu(0))\|_{H_0^1(\Omega)\times L^2(\Omega)}
\]
 In addition, solution $u$ to \eqref{NLW} satisfies an energy conservation law \cite{TT}:
 \[
 E(u)(t)=\int_\Omega\left(\frac{1}{2}|\partial_t u|^2+\frac{1}{2}|\nabla u|^2+\frac{1}{6}|u|^6\right)dx=E(u)(0).
 \]
 We remark that this explains why the exponent $5$ in the nonlinear term $u^5$ is the $H^1$ critical wave since it satisfies $5=\frac{d+2}{d-2}$ with the space dimension $d=3$.
By the Sobolev embedding $H_0^1(\Omega)\hookrightarrow L^6(\Omega)$,  we see that the energy $E(u)(0)$ is finite for any initial data $(u_0,u_1)\in H_0^1(\Omega)\times L^2(\Omega).$

The global existence  for smooth solutions for energy critical wave \eqref{NLW} was proven for $\Omega=\mathbb{R}^3$ in \cite{MG} and for energy space solutions in \cite{ShS,ShS1}. This result is extended to the case of the exterior of convex obstacles in \cite{SS2}. Finally, the case of bounded domains $\Omega\subset\mathbb{R}^3$ with Dirichlet conditions was treated in \cite{BLP}. 

Let us also recall some results on dispersive and Strichartz estimates as they have been proved to powerful in the study of well-possedness of nonlinear problems. 

The dispersive estimates for the wave equation in $\mathbb{R}^ d$  follows from the representation of the solution as a sum of Fourier integral operators (see \cite{BCD,PB,GV}). They read as follows:
\begin{align}\label{freespace}
\lVert\chi(hD_t)e^{\pm it\sqrt{-\Delta_{\mathbb{R}^d}}}\rVert_{L^1(\mathbb{R}^ d)\rightarrow L^\infty(\mathbb{R}^ d)}\leq Ch^{-d}\min\Bigg\{1,\bigg(\frac{h}{|t|}\bigg)^{\frac{d-1}{2}}\Bigg\},
\end{align}
where $\Delta_{\mathbb{R}^d}$ is the Laplace operator in $\mathbb{R}^d$. Here and in the sequel, the function $\chi$ belongs to $C_0^\infty(]0,\infty[)$ and is equal to $1$ on $[1,2]$ and $D_t=\frac{1}{i}\partial_t$. 

Inside strictly convex domains $\Omega_D$ of dimensions $d\geq 2$, the optimal (local in time) dispersive estimates for the wave equations have been established  by Ivanovici, Lebeau, and Planchon in \cite{ILP}. More precisely, they have proved that 
 \begin{align}\label{domain}
\lVert\chi(hD_t)e^{\pm it\sqrt{-\Delta_D}}\rVert_{L^1(\Omega_D)\rightarrow L^\infty(\Omega_D)}\leq Ch^{-d}\min\Bigg\{1,\bigg(\frac{h}{|t|}\bigg)^{\frac{d-1}{2}-\frac{1}{4}}\Bigg\},
\end{align}
 where $\Delta_D$ is the Laplace operator on $\Omega_D$.
  Due to the caustics formation in arbitrarily small times, \eqref{domain} induces a loss of $\frac{1}{4}$ powers of $(h/|t|)$ factor compared to the free wave estimates \eqref{freespace}.
  
The Strichartz estimates on Riemannian manifold (see \cite{ILP}, section 1).
Let $(\Omega,g)$ be a Riemannian manifold  without boundary of dimensions $d\geq 2$. Local in time Strichartz estimates state that 
 \begin{align}\label{stri}
 \lVert u\rVert_{L^q((-T,T);L^r(\Omega))}\leq C_T\bigg(\lVert u_0\rVert_{\dot{H}^\beta(\Omega)}+\lVert u_1\rVert_{\dot{H}^{\beta-1}(\Omega)}\bigg),
 \end{align}
where $\dot{H}^\beta$ is the  homogeneous Sobolev space over $\Omega$ of order $\beta$ and $2\leq q,r\leq \infty$ satisfy
\[
\frac{1}{q}+\frac{d}{r}=\frac{d}{2}-\beta,\quad \frac{1}{q}\leq \frac{d-1}{2}\bigg(\frac{1}{2}-\frac{1}{r}\bigg).
\]
Here $u=u(t,x)$ is a solution to the wave equation
\begin{align*}
(\partial_t^2-\Delta_g)u=0\,\,\text{in } (-T,T)\times \Omega, \quad u(0,x)=u_0(x),\quad \partial_tu(0,x)=u_1(x),
\end{align*} 
where $\Delta_g$ denotes the Laplace-Beltrami operator on $(\Omega,g)$. The estimates \eqref{stri} hold on $\Omega=\mathbb{R}^d$ and $g_{ij}=\delta_{ij}.$

 In \cite{BSS2}, Blair, Smith, Sogge proved the Strichartz estimates for the wave equation on (compact or noncompact) Riemannian manifold with boundary. They proved that the Strichartz estimates \eqref{stri} hold if $\Omega$ is a compact manifold with boundary
and $ (q,r,\beta)$ is a triple satisfying
\[
\frac{1}{q}+\frac{d}{r}=\frac{d}{2}-\beta \quad\text{for}\quad 
\begin{cases}
\frac{3}{q}+\frac{d-1}{r}\leq\frac{d-1}{2},&d\leq 4,\\
\frac{1}{q}+\frac{1}{r}\leq \frac{1}{2},&d\geq 4.
\end{cases}
\]

Recently in \cite{ILP}, Ivanovici, Lebeau, and Planchon have deduced a local in time Strichartz estimates \eqref{stri} from the optimal dispersive estimates inside strictly convex domains of dimensions $d\geq 2$
for a triple $(d,q,\beta)$ satisfying
\[
\frac{1}{q}\leq\bigg(\frac{d-1}{2}-\frac{1}{4}\bigg)\bigg(\frac{1}{2}-\frac{1}{r}\bigg)
\,\,\text{ and }\,\, \beta=d\bigg(\frac{1}{2}-\frac{1}{r}\bigg)-\frac{1}{q}.
\]
 For $d\geq 3$ this improves the range of indices for which sharp Strichartz estimates do hold compared to the result by   Blair, Smith, Sogge in \cite{BSS2}. However, the results in \cite{BSS2} apply to any domains or manifolds with boundary.
 
The latest results in \cite{ILP4} on Strichartz estimates inside the Friedlander model domain have been obtained for pairs $(q, r)$ such that 
 $$\frac{1}{q}\leq \bigg(\frac{1}{2}-\frac{1}{9}\bigg)\bigg(\frac{1}{2}-\frac{1}{r}\bigg).$$
This result improves on the known results for strictly convex domains for $d=2$, while in \cite{ILP} only gives a loss of $\frac{1}{4}$.

In our case, we treat the equation \eqref{NLW} in $\Omega\subset\mathbb{R}^3$, where $\Omega$ is locally a cylindrical convex domains. The main result of this work on the local well-posedness is stated as the following.

\begin{theorem}\label{LOC}
Let $(\Omega, \Delta)$ be defined as before. For any $(u_0, u_1)\in H_0^1(\Omega)\times L^2(\Omega)$ there exists a local solution $u$ to \eqref{NLW} in the space
 \[
 C^0((0,T); H_0^1(\Omega))\cap C^1((0,T); L^2(\Omega))\cap L^5((0,T); L^{10}(\Omega)).
 \]
\end{theorem}
We note that the nonlinearity is defocusing due to its sign, which does not play a role for the local existence of solutions to \eqref{NLW} and hence the result in Theorem \ref{LOC} also hold in case of focusing quintic wave equation. While the sign of nonlinear term is crucial for global in time existence of solutions to \eqref{NLW}. The main difficulty for proving global in time solutions to \eqref{NLW} is that one does not obtain a bound on $u\in L^5((0,T);L^{10}(\Omega))$ and thus on $u^5\in L^1((0,T);L^2(\Omega))$. But in our domain, the Strichartz estimates in Theorem \ref{STR} allows us to control $L^5((0,T);L^{10}(\Omega))$ of the solution to  \eqref{NLW} by the energy norm of the initial data.

The Theorem \ref{LOC} was established  by Burq, Lebeau, Planchon in  \cite{BLP} for any bounded domains $\Omega_D\subset\mathbb{R}^3$. Their idea is based on the spectral projector established by Smith and Sogge in \cite{SS3} to derive the optimal and scale invariant Strichartz estimates for the solution to the wave equation, while in our setting we will use the dispersive estimates obtained in \cite{L,L1,L3} and the Strichartz estimates in \cite{L2} to establish local existence of solutions to \eqref{NLW}.  

We point out that the approach in \cite{BLP} to get the Strichartz estimates to control the $L^5\big((0,1); W_0^{\frac{3}{10},5}(\Omega_D)\big)$ norm of the solution of the wave equation by the energy norm, assuming $\Omega$ is compact. This estimate allows them to extend local to global well-posedness for arbitrary finite energy initial data, when combined with the estimate for the normal derivative and the nonconcentration estimate of nonlinear effect. 

The global well-posedness in our setting reads as follows.

\begin{theorem}\label{GWP}
Let $(\Omega, \Delta)$ be defined as before. For any $(u_0, u_1)\in H_0^1(\Omega)\times L^2(\Omega)$ there exists a unique global solution $u$ to \eqref{NLW} in the space
\[
C^0(\mathbb{R}_t; H_0^1(\Omega))\cap C^1(\mathbb{R}_t; L^2(\Omega))\cap L_{loc}^5(\mathbb{R}_t; L^{10}(\Omega)).
\]
\end{theorem}

In this paper, for $s\geq 0$,  let $\dot{H}^{s}(\Omega)=(-\Delta )^{-\frac{s}{2}}L^2(\Omega)$ be the homogeneous Sobolev space over $\Omega$ and $H_0^s(\Omega)$ is the closure in $\dot{H}^{s}(\Omega)$ of the set of smooth and compactly supported functions. We note that $\dot{H}^{s}(\Omega)=H_0^s(\Omega)$ for $0\leq s<\frac{3}{2}$, and that when $s=1, H_0^1(\Omega)$ is a Hilbert space with the inner product of $H^1(\Omega)$. The notation $A\lesssim B$ means that there exists a constant $C$ such that $A\leq CB$ and this constant may change from line to line but is independent of all parameters. Similarly, $A\sim B$ means there exist constants $C_1, C_2$ such that $C_1B\leq A\leq C_2B$.

\section{Dispersive Estimates}
In this section, we will present the result on the local in time dispersive estimates for the solution to the linear wave equation in the cylindrical domain $\Omega$ with the Laplace $\Delta$ defined as before and we obtained a sharp loss of $\frac{1}{4}$ due to swallowtail type singularities in the wave front set.

Let consider the Dirichlet wave equation inside $\Omega=\{x\geq0,(y,z)\in\mathbb{R}^2\}\subset\mathbb{R}^3$
\begin{equation}\label{WD}
\begin{aligned}
&(\partial_t^2-\Delta)u=0\quad \text{in } \mathbb{R}_t\times\Omega,\\
&u_{|t=0}=\delta_a,\,\,\, \partial_t u_{|t=0}=0,\,\,\, u_{|x=0}=0.
\end{aligned}
\end{equation}
where the Dirac distribution $\delta_a=\delta_{x=a,y=0,z=0}$ with $(a,0,0)\in\Omega, a>0$. In local coordinates $a$ denotes the distance from the source point to the boundary of $\Omega$. We assume that $0<a\ll 1$ is small enough as we are interested only in highly reflected waves, which give us interesting phenomena such as caustics near the boundary.

The local in time dispersive estimates for the wave equation \eqref{WD} are established in \cite{L,L1,L3}. Let $\chi\in C_0^\infty(]0,\infty[), \chi=1$ on $[1,2]$ and $D_t=\frac{1}{i}\partial_t$. Let $\mathcal{G}_a$ be the Green function for \eqref{WD}.
\begin{theorem}\label{DISP}
There exists $C$ such that for every $h\in]0,1]$, every $t\in [-1,1]$ and every $a\in ]0,1]$ the following holds:
\[
\|\chi(hD_t)\mathcal{G}_{a}(t,x,y,z)\|_{L^1(\Omega)\to L^\infty(\Omega)}\leq C h^{-3}\min\bigg\{1, \bigg(\frac{h}{|t|}\bigg)^{\frac{3}{4}}\bigg\}.
\]
\end{theorem}
In our cylindrical domain with boundary, the light rays may no longer slightly distorted straight lines. There may be rays glancing near tangential direction of the boundary or rays gliding along a convex part of the cylinder boundary or combinations of both.  We
note that the interesting phenomena analyzed in \cite{L,L1,L3} is the caustics (cusps and swallowtails) near the boundary and due to these caustics on the boundary, the dispersive estimates in Theorem \ref{DISP} has a sharp loss of $\frac{1}{4}$ powers of $(h/|t|)$ factor compared to the free wave estimates in dimension $3$. This is compatible with the intuition: near the boundary less dispersion occurs compared to the $\mathbb{R}^3$ case. Moreover, the geometry analysis of the wave front set allows us to track the swallowtail type singularities in the Green function originating at $x=a$ in the $(x,t)$ plane. The singular points are $(a, 4N\sqrt{a(1+a)/(1-\delta^2)})$ for $|N|\lesssim 1/\sqrt{a}, |\delta|<1,$ where the times depend on the frequency of the source and its distance to the boundary. Let us mention that if $\delta=0$, it is exactly where the swallowtail singularity for the Friedlander's model in \cite{ILP} occurs.

Let us sketch the proof and key ingredients of Theorem \ref{DISP} of the frequency-localized dispersive estimates as follows (see \cite{L3}): the strategy is based on the construction of parametrices for the fundamental solution of the wave equation \eqref{WD} and possibly degenerate stationary phase method.

A key feature of our Laplacian $\Delta=\partial_x^2+(1+x)\partial_y^2+\partial_z^2$ on the half space $\Omega$ is that the coefficients of the metric do not depends on the variables $y$ and $z$ and therefore this allows us to take the Fourier transform in $y$ and $z$. Now taking the Fourier transform in $y,z$-variables yields
\[
-\Delta_{\eta,\zeta}:=-\partial_x^2+(1+x)\eta^2+\zeta^2.
\]

For $\eta\neq0,-\Delta_{\eta,\zeta}$ is a self-adjoint, positive operator on $L^2(\mathbb{R}_+)$ with a compact resolvent.
Let $(e_k)_{k\geq 1}$ be an orthonormal basis in $L^2(\mathbb{R}_+)$ of Dirichlet eigenfunctions  of $ -\Delta_{\eta,\zeta}$ and let $(\lambda_k)_k$ be the associated eigenvalues. That is,
\[
-\Delta_{\eta,\zeta}e_k=\lambda_k e_k.
\]
These eigenfunctions are explicit in term of Airy functions
  \begin{align*}
  e_k=&e_k(x,\eta)=f_k\frac{|\eta|^{1/3}}{k^{1/6}}Ai(|\eta|^{2/3}x-\omega_k)
\end{align*}
with associated eigenvalues
\begin{align*}
\lambda_k=&\lambda_k(\eta,\zeta)=\eta^2+\zeta^2+\omega_k|\eta|^{4/3},
\end{align*}
where $(-\omega_k)_k$ denote the zeros of Airy function in decreasing order and for all $k\geq 1,$ $f_k$ are constants so that $\lVert e_k(.,\eta)\rVert_{L^2(\mathbb{R}_+)}=1$. Observe that $(f_k)_k$ is uniformly bounded in a fixed compact subset of $]0,\infty[$ as a consequence of  \[\int_{-\omega_k}^{-2}Ai^2(\omega)d\omega\sim\frac{1}{4\pi}\int_{-\omega_k}^{-2}|\omega|^{-1/2}(1+O(\omega^{-1}))d\omega\sim |\omega_k|^{1/2}\] and \[\omega_k\sim\bigg(\frac{3}{2}\pi k\bigg)^{2/3}(1+O(k^{-1})).\]

 For $a\in\Omega$, let $g_a(t,x,\eta,\zeta)$ be the solution of
\begin{align*}
(\partial_t^2+\Delta_{\eta,\zeta})g_a=0,\,\,
{g_a}_{|x=0}=0,\,\, {g_a}_{|t=0}=\delta_{x=a},\,\,\partial_t{g_a}_{|t=0}=0.
\end{align*}
We have
\begin{align}
g_a(t,x,\eta,\zeta)=\sum_{k\geq1}\cos(t\lambda_k^{1/2})e_k(x,\eta)e_k(a,\eta).
\end{align}
Here $\delta_{x=a}$ denotes the Dirac distribution on $\mathbb{R}_+$ with $a>0$ and it may be decomposed in terms of eigenfunctions $(e_k)_{k\geq 1}$ as follows:
 \[\delta_{x=a}=\sum\limits_{k\geq1}e_k(x,\eta)e_k(a,\eta).\] 
Now applying the inverse Fourier transform, the Green function for \eqref{WD} is given by
\begin{align}
\mathcal{G}_a(t,x,y,z)&=\frac{1}{4\pi^2}\int e^{i(y\eta+z\zeta)}g_a(t,x,\eta,\zeta)d\eta d\zeta,\nonumber\\
&=\frac{1}{4\pi^2h^2}\sum_{k\geq1}\int e^{i(y\eta+z\zeta)/h}\cos(t\lambda_k^{1/2})e_k(x,\eta/h)e_k(a,\eta/h)d\eta d\zeta.
\end{align}
We perform a spectral localization $\lambda_k\sim h^{-2}$, corresponding to inserting a smooth, compactly supported away from zero $\chi(h\sqrt{\lambda_k})$; on the wave flow, this is $\chi(hD_t)$ and it smooths out the Green function. We thus get the formula for $2\chi(hD_t)\mathcal{G}_a$ as follows.
\begin{align}
2\chi(hD_t)\mathcal{G}_a(t,x,y,z)&=\frac{1}{4\pi^2h^2}\sum_{k\geq1}\int 
e^{\frac{i}{h}(y\eta+z\zeta)}e^{i\frac{t}{h}(\eta^2+\zeta^2+\omega_kh^{2/3}|\eta|^{4/3})^{1/2}} e_k(x,\eta/h)\nonumber\\& \quad \times e_k(a,\eta/h)\chi((\eta^2+\zeta^2+\omega_kh^{2/3}|\eta|^{4/3})^{1/2})d\eta d\zeta.
\end{align}
On the wave front set of the above expression, one has $\tau=(\eta^2+\zeta^2+\omega_kh^{2/3}|\eta|^{4/3})^{1/2}$. 
In order to prove Theorem \ref{DISP}, we only need to work near tangential directions; therefore, we will introduce an extra cutoff
to insure $\vert \tau-(\eta^2+\zeta^2)^{1/2}\vert$ small, which is equivalent to  $\omega_kh^{2/3}|\eta|^{4/3}$ small.
Therefore, we are reduced to prove the dispersive estimate for $\mathcal{G}_{a,loc}$:

\begin{align}\label{eq:kparametrix}
\mathcal{G}_{a,loc}(t,x,y,z)&=\frac{1}{4\pi^2h^2}\sum_{k\geq1}\int 
e^{\frac{i}{h}(y\eta+z\zeta)}e^{i\frac{t}{h}(\eta^2+\zeta^2+\omega_kh^{2/3}|\eta|^{4/3})^{1/2}} e_k(x,\eta/h)\nonumber\\&\quad \times e_k(a,\eta/h)\chi_0(\eta^2+\zeta^2)\chi_1(\omega_kh^{2/3}|\eta|^{4/3}) d\eta d\zeta,
\end{align}
where the cutoff functions $\chi_0\in C_0^\infty, 0\leq\chi_0\leq 1, \chi_0$ is supported in the neighborhood of $1$ and $\chi_1\in C_0^\infty, 0\leq\chi_1\leq 1, \chi_1$ is supported in $(-\infty, 2\varepsilon],\chi_1=1$ on $(-\infty,\varepsilon].$

We remark that \eqref{eq:kparametrix} is a parametrix in the frequency localization and near tangential direction which is the sum of oscillatory integrals with phase functions containing  Airy type functions with degenerate critical points. We give a precise analysis of the Lagrangian in the phase space associated to these oscillatory integrals. This geometric analysis allows us to track the degeneracy of the phases when we apply the stationary phase method.

We further analyze \eqref{eq:kparametrix} to get the local in time dispersive estimates by cutting the $\eta$ integration into different pieces. More precisely, we write
\[
\mathcal{G}_{a,loc}=\mathcal{G}_{a,c_0}+\sum_{\epsilon_0\leq 2^m\sqrt{a}\leq c_0}\mathcal{G}_{a,m}+\mathcal{G}_{a,\epsilon_0},
\] 
where $\mathcal{G}_{a, c_0}$ is associated with the integration for $|\eta|\geq c_0$, $\mathcal{G}_{a,m}$ is associated with the integration for $|\eta|\sim 2^m\sqrt{a}$ and $\mathcal{G}_{a,\epsilon_0}$ is associated with the integration for $0<|\eta|\leq\epsilon_0\sqrt{a}$.

We have following results established in \cite{L1,L3}. Let $\epsilon\in ]0,1/7[$.

\begin{theorem}\label{betabis}
There exists $C$ such that for every $h\in ]0,1]$, every $t\in [h,1]$,  the following holds:
\begin{align}
\lVert\mathcal{G}_{a,c_0}(t,x,y,z)\rVert_{L^\infty(x\leq a)}\leq Ch^{-3}\bigg(\frac{h}{t}\bigg)^{1/2}\gamma(t,h,a),
\end{align}
with 
\[
\gamma(t,h,a)=\begin{cases}\Big(\frac{h}{t}\Big)^{1/3}&\text{if}\,\,\, a\leq h^{\frac{2}{3}(1-\epsilon)},\\
\Big(\frac{h}{t}\Big)^{1/2}+a^{1/8}h^{1/4}&\text{if}\,\,\, a\geq h^{\frac{2}{3}(1-\epsilon')},\epsilon'\in ]0,\epsilon[
.\end{cases}
\]
\end{theorem}
The key ingredient is to use the local parametrix $\mathcal{G}_{a,c_0}$, which is a sum of the eigenmodes (over $k$) to prove the estimates for $a\leq h^{\frac{2}{3}(1-\epsilon)},\epsilon\in]0,1/7[$. Using the Airy-Poisson summation formula (see Lemma 2.4 \cite{L3}), $\mathcal{G}_{a, c_0}$ can be also represented as a sum over $N\in\mathbb{Z}$ for $a\geq h^{\frac{2}{3}(1-\epsilon')},$ for $\epsilon'\in]0,\epsilon[$. These are oscillatory integrals with the phase functions containing Airy type functions with degenerate critical points and to establish their estimates we use the stationary phase method (possibly degenerate) together with the geometric analysis of the Lagrangian associated to the phase functions. We note that on the support of $\chi_1$, we have $\omega_k h^{2/3}|\eta|^{4/3}\leq 2\varepsilon$ and since $\omega_k\sim k^{2/3}$; we obtain $k\leq \frac{\varepsilon}{h|\eta|^2}$. Thus since $\eta$ is bounded from below, we may assume that $k\leq\varepsilon/h.$ Similarly, the geometric study of the wave front set and the Lagrangian submanifold corresponding to the phase function of $\mathcal{G}_{a,c_0}$ when it is expressed as the sum over $N\in\mathbb{Z}$ allows us to reduce the sum over $N\in\mathbb{Z}$ to the sum over $1\leq N\leq Ca^{-1/2}$.

\begin{theorem}\label{1beta}
There exists $C$ such that for every $h\in ]0,1]$, every $t\in [h,1]$, the following holds:
\begin{equation*}
\lVert\mathcal{G}_{a,m}(t,x,y,z)\rVert_{L^\infty(x\leq a)}\leq Ch^{-3}\left(\frac{h}{t}\right)^{1/2}\gamma_m\big(t,h,a),
\end{equation*}
with 
\[
\gamma_m(t,h,a)=\begin{cases}
\Big(\frac{h}{t}\Big)^{1/3}(2^m\sqrt a)^{1/3}\qquad\qquad\quad\text{ if } a\leq \bigg(\frac{h}{2^m\sqrt a}\bigg)^{\frac{2}{3}(1-\epsilon)},\\
\min \bigg\{\Big(\frac{h}{t}\Big)^{1/2}, 2^m\sqrt a\vert \log(2^m\sqrt a)\vert \bigg\}+a^{1/8}h^{1/4}(2^m\sqrt a)^{3/4}\\
\hspace{4.5cm}\text{ if } a\geq \!\bigg(\frac{h}{2^m\sqrt a}\bigg)^{\frac{2}{3}(1-\epsilon')}\!\!\!,
\epsilon'\in ]0,\epsilon[.\\\end{cases}
\]
\end{theorem}
 For $2^m\sqrt a\sim 1$, Theorem \ref{1beta} yields the same result as in Theorem \ref{betabis}. We notice that the estimates get better when $|\eta|$( $\sim 2^m\sqrt a)$ decreases. This is compatible with the intuition that less curvature implies better dispersion. 
 
 In this regime of $\eta$, we apply the same approach as in the case of $\eta$ is bounded from below. But we have to handle carefully in the case when $2^m\sqrt{a}$ could be very small. As proved in Lemma 3.3 \cite{L3}, whenever $2^m\sqrt{a}\lesssim h$, there is another factor $2^m\sqrt{a}|\log(2^m\sqrt{a})|$ in the estimates.

 \begin{theorem}\label{0eta}
There exists $C$ such that for every $h\in ]0,1]$, every $t\in [h,1],$  the following holds:
\begin{align}
\lVert\mathcal{G}_{a,\epsilon_0}(t,x,y,z)\rVert_{L^\infty(x\leq a)}\leq Ch^{-3} \bigg(\frac{h}{t}\bigg)^{1/2} \min \bigg\{ \bigg(\frac{h}{t}\bigg)^{1/2} , \sqrt a\vert \log(a)\vert  \bigg\}.
\end{align}
\end{theorem}
It is particularly interesting at this localization $\mathcal{G}_{a,\epsilon_0}$ is an oscillatory integral with nondegenerate phase function; this is due to the geometric study of the associated Lagrangian of the characteristic set of the wave operator which rules out the swallowtails singularity for $|t|<1$ if $\epsilon_0$ is mall enough. Then we follow the preceding technique to estimate the oscillatory integrals. But we point out that a new key ingredient in this localization; that is, the Lemma 4.3 in \cite{L3} which necessary in the geometric study of the wave front set and the Lagrangian submanifold associated to the phase function.  

 Let us verify that Theorem \ref{DISP} is a consequence of Theorems \ref{betabis}, \ref{1beta} and \ref{0eta}.
  We may assume $\vert t\vert\geq h$, since for $\vert t\vert\leq h$, the best bound for the dispersive estimate is equal to $Ch^{-3}$ by Sobolev inequality. Then, by symmetry of the Green function, we may assume $t\in [h,1]$
and $x\leq a$. Then Theorem \ref{DISP} is a consequence of $\sum_{m\leq M}(2^m\sqrt a)^\nu\sim (2^M\sqrt a)^\nu$
for $\nu>0$.

We remark that on the wave flow Theorem \ref{DISP} holds true if we replace $\chi(hD_t)\mathcal{G}_a$ by $\chi(h\sqrt{-\Delta})\mathcal{G}_{a}$. This is useful as it allows us to connect the definition of the homogeneous Sobolev space norm defined by Littlewood-Paley decomposition and frequency-localized dispersive estimates.

\section{Strichartz Estimates}
In this section, we show the Strichartz estimates \label{STR} inside cylindrical convex domains in dimension $3$. The technique follows from the abstract Strichartz estimates proved in Keel-Tao \cite{KT}. Let us recall the semiclassical Strichartz estimates version (see \cite{CDS,ZJ,MZ})
\begin{proposition}[Semiclassical Strichartz Estimates]
Let $(X,\mathcal{M}, dm)$ be a $\sigma$-finite measured space and $U\in L^\infty(\mathbb{R};\mathcal{B}(L^2(X,\mathcal{M},dm))$ satisfies
\begin{align*}
&\|U(t)\|_{L^2\to L^2}\leq C,\,\,t\in\mathbb{R},\\
&\|U(t)U^*(s) f\|_{L^\infty}\leq C h^{-\alpha}|t-s|^{-\sigma}\|f\|_{L^1},
\end{align*}
some constants $C,\alpha\geq 0,\sigma>0, h>0$.
The for every pair $q,r\in [1,\infty]$ such that $(q,r,\sigma)\neq (2,\infty,1)$ and
\[
\frac{1}{q}+\frac{\sigma}{r}\leq\frac{\sigma}{2},\,\, q>2,
\]
there exists a constant $\tilde C$ only depending on $C,\sigma,q,$ and $r$ such that
\[
\left(\int_{\mathbb{R}}\|U(t)f\|_{L^r}^qdt\right)^{\frac{1}{q}}\leq \tilde Ch^{-(\alpha+\sigma)\big(\frac{1}{2}-\frac{1}{r}\big)+\frac{1}{q}}\|f\|_{L^2}.
\]
\end{proposition}

The frequency-localized dispersive estimates in Theorem \ref{DISP} yields the frequency-localized Strichartz estimates for the wave equation inside cylindrical convex domains in dimension $3$ by conservation of energy, interpolation and $TT^*$ arguments.

Let us recall the Littlewood-Paley decomposition and some links with Sobolev spaces \cite{BCD}.
Let $\chi\in C_0^\infty(\mathbb{R}^*)$ and equal to $1$ on $[1/2,2]$ such that 
\[
\sum_{j\in\mathbb{Z}}\chi(2^{-j}\lambda)=1,\,\,\, \lambda>0.
\]
We define the associated Littlewood-Paley frequency cutoffs $\chi(2^{-j}\sqrt{-\Delta})$ using spectral theorem for $\Delta$ and we have 
\[
\sum_{j\in\mathbb{Z}}\chi(2^{-j}\sqrt{-\Delta})=\text{Id}: L^2(\Omega)\longrightarrow L^2(\Omega).
\]

This decomposition takes a single function and writes it as a superposition of a countably infinite family of functions $\chi$ each one having a frequency of magnitude $\sim 2^{j}$, for $j\geq 1$. We have that a norm of the homogeneous Sobolev of $\dot{H}^{\beta}$ is defined as follows:
for all $\beta\geq 0$,
\begin{align*}
\|u\|_{\dot{H}^\beta(\Omega)}:=\left(\sum_{j\in\mathbb{Z}}2^{2j\beta}\|\chi(2^{-j}\sqrt{-\Delta})u\|_{L^2(\Omega)}^2\right)^{1/2}.
\end{align*}
With this decomposition, the result about the Littlewood-Paley squarefunction estimate (see \cite{BFHM,BFM,OF}) reads as follows:
for $f\in L^r(\Omega),\forall r\in [2,\infty[$,
\begin{align}\label{sqrf}
\|f\|_{L^r(\Omega)}\leq C_r\left\|\left(\sum_{j\in\mathbb{Z}}|\chi(2^{-j}\sqrt{-\Delta})f|^2\right)^{1/2}\right\|_{L^r(\Omega)}.
\end{align}

We define the frequency localization $v_j$ of $u$ by $v_j=\chi(2^{-j}\sqrt{-\Delta})u$. Hence, $u=\sum_{j\geq 0}v_j$.
 Let $h=2^{-j}$. We deduce from the dispersive estimates inside cylindrical convex domains established in Theorem \ref{DISP} the frequency-localized dispersive estimates as follows.
 \begin{proposition}
 Let  $v_j=\chi(hD_t)u$ be the solution to the (frequency-localized) wave equation
\begin{equation}\label{eq:fleq}
\begin{aligned}
&(\partial_t^2-\Delta)v_j=0\,\,\, \text{in}\,\,\, \mathbb{R}_t\times\Omega,\\
&{v_j}_{|t=0}=\chi(hD_t)u_0,\,\,
\partial_t{v_j}_{|t=0}=\chi(hD_t)u_1,\,\,
{v_j}_{|\partial\Omega}=0.
\end{aligned}
\end{equation}
Then it holds true,
\begin{align}
&\left\|\cos(t\sqrt{-\Delta})\chi(hD_t)u_0\right\|_{L^\infty(\Omega)}\lesssim h^{-3}\min\left\{1,\Big(\frac{h}{t}\Big)^{\frac{3}{4}}\right\}\|\chi(hD_t)u_0\|_{L^1(\Omega)}\label{cos}\\
&\left\|\frac{\sin(t\sqrt{-\Delta})}{\sqrt{-\Delta}}\chi(hD_t)u_1\right\|_{L^\infty(\Omega)}\lesssim h^{-2}\min\left\{1,\Big(\frac{h}{t}\Big)^{\frac{3}{4}}\right\}\|\chi(hD_t)u_1\|_{L^1(\Omega)}\label{sin}
\end{align}
and the energy estimates
\begin{align}
&\left\|\cos(t\sqrt{-\Delta})\chi(hD_t)u_0\right\|_{L^2(\Omega)}\lesssim\|\chi(hD_t)u_0\|_{L^2(\Omega)}\label{cos1}\\
&\left\|\frac{\sin(t\sqrt{-\Delta})}{\sqrt{-\Delta}}\chi(hD_t)u_1\right\|_{L^2(\Omega)}\lesssim\|\chi(hD_t)u_1\|_{L^2(\Omega)}.\label{sin1}
\end{align}
In particular, interpolating \eqref{cos} with \eqref{cos1} (respectively \eqref{sin} with \eqref{sin1}) yields
\begin{align*}
&\left\|\cos(t\sqrt{-\Delta})\chi(hD_t)u_0\right\|_{L^r(\Omega)}\lesssim h^{-3\left(1-\frac{2}{r}\right)}\min\left\{1,\Big(\frac{h}{t}\Big)^{\frac{3}{4}\left(1-\frac{2}{r}\right)}\right\}\|\chi(hD_t)u_0\|_{L^{r'}(\Omega)}\\
&\left\|\frac{\sin(t\sqrt{-\Delta})}{\sqrt{-\Delta}}\chi(hD_t)u_1\right\|_{L^r(\Omega)}\lesssim h^{-2\left(1-\frac{2}{r}\right)}\min\left\{1,\Big(\frac{h}{t}\Big)^{\frac{3}{4}\left(1-\frac{2}{r}\right)}\right\}\|\chi(hD_t)u_1\|_{L^{r'}(\Omega)},
\end{align*}
for all $r\in[2,\infty]$, where $r'$ denotes the exponent conjugate to $r$; that is, $\frac{1}{r}+\frac{1}{r'}=1$.
\end{proposition}
The $L^1\to L^\infty$ estimates, the energy estimates  for the wave equation, the Riesz-Thorin interpolation, and $TT^*$ argument yield the following frequency-localized Strichartz estimates proved in \cite{L2}.
\begin{theorem}[Frequency-Localized Strichartz Estimates]
Let $(\Omega,\Delta)$ be defined as before. 
Let $v_j$ be a solution of the (frequency-localized) wave equation \eqref{eq:fleq}.
Then it holds true,
\begin{align}
h^\beta\left\|\cos(t\sqrt{-\Delta})\chi(hD_t)u_0\right\|_{L^q((0,T);L^r(\Omega))}\lesssim \|\chi(hD_t)u_0\|_{L^2(\Omega)},\label{ssin}\\
h^{\beta-1}\left\|\frac{\sin(t\sqrt{-\Delta})}{\sqrt{-\Delta}}\chi(hD_t)u_1\right\|_{L^q((0,T);L^r(\Omega))}\lesssim \|\chi(hD_t)u_1\|_{L^2(\Omega)},\label{scos}
\end{align}
with \[q\in ]2,\infty],\,\,r\in[2,\infty],\,\,\frac{1}{q}\leq\frac{3}{4}\left(\frac{1}{2}-\frac{1}{r}\right),
\,\,\,\text{and}\,\,\, \beta=3\left(\frac{1}{2}-\frac{1}{r}\right)-\frac{1}{q}.\]
\end{theorem}
Let us make it precise for \eqref{ssin} as \eqref{scos} we may treat analogously. By $TT^*$ argument, one has \eqref{ssin} is equivalent to
\[
h^\beta\left\|\int_0^T \cos(t\sqrt{-\Delta})\chi(hD_t)F(t,.)dt\right\|_{L^2(\Omega)}\lesssim\|\chi(hD_t)F\|_{L^{q'}((0,T);L^{r'}(\Omega))}
\]
and both are equivalent to
\[
h^{2\beta}\left\|\int_0^T \cos((t-s)\sqrt{-\Delta})\chi(hD_t)F(s,.)ds\right\|_{L^q((0,T);L^r(\Omega))}\lesssim\|\chi(hD_t)F\|_{L^{q'}((0,T);L^{r'}(\Omega))}.
\]
 The establishment of the Strichartz estimates in the context of cylindrical domains requires the Littlewood-Paley square function estimates (see \cite{OF}). The following result is established in \cite{L2}.
\begin{theorem}[Homogeneous Strichartz Estimates]\label{STR}
Let $(\Omega,\Delta)$ be defined as before. Let $u$ be a solution of the following wave equation on  $\Omega$:
\begin{equation}\label{STRW}
\begin{aligned}
&(\partial_t^2-\Delta)u=0\,\,\, \text{in}\,\,\,\mathbb{R}_t\times \Omega,\\
&u_{|t=0}=u_0,\,\,\, \partial_t u_{|t=0}=u_1,\,\,\,u_{|x=0}=0,
\end{aligned}
\end{equation}
for some initial data $(u_0,u_1)\in \dot{H}^\beta(\Omega)\times\dot {H}^{\beta-1}(\Omega)$.
 Then
\[
\|u\|_{L^q((0,T);L^r(\Omega))}\leq C_{\text{Str}}\left(\|u_0\|_{\dot{H}^{\beta}(\Omega)}+\|u_1\|_{\dot{H}^{\beta-1}(\Omega)}\right),
\]
with 
\begin{align}\label{ADM}
\frac{1}{q}\leq\frac{3}{4}\left(\frac{1}{2}-\frac{1}{r}\right)
\end{align}
and
\begin{align}\label{ADM1}
 \beta=3 \left(\frac{1}{2}-\frac{1}{r}\right)-\frac{1}{q}.
\end{align}
\end{theorem}
Remark that the Strichartz estimates in Theorem \ref{STR} in dimension $3$ improves the range of indices for which the sharp Strichartz estimates hold compared to the result in \cite{BSS2}. However, the result in Theorem \ref{STR} is restricted to cylindrical domains, while \cite{BSS2} applies to any domains or manifolds with boundary.

We note also that from the result in Theorem \ref{STR}, we have control of the $L^4((0,T); L^{12}(\Omega))$ norms and the  $L^5((0,T); L^{10}(\Omega))$ of the solution of the wave equation in terms of the energy norm of the initial data, which are important estimates to prove that there is scattering for \eqref{NLW} when the domain is the compliment of a star-shaped obstacle \cite{BSS2}.

Finally, the inhomogeneous Strichartz estimates follow from the homogeneous Strichartz estimates and the Christ-Kiselev lemma \cite{CK}. In the following corollary $(\tilde q', \tilde r')$ denotes the exponents conjugate to $(\tilde q, \tilde r)$.
\begin{corollary}\label{INS}
Let $(\Omega, \Delta)$ be defined as before. Let $u$ be a solution of the following wave equation on  $\Omega$:
\begin{equation}
\begin{aligned}
&(\partial_t^2-\Delta)u=F\,\,\, \text{in}\,\,\,\mathbb{R}_t\times \Omega,\\
&u_{|t=0}=u_0,\,\,\, \partial_t u_{|t=0}=u_1,\,\,\,u_{|x=0}=0.
\end{aligned}
\end{equation}
for some initial data $(u_0,u_1)\in\dot{H}^\beta(\Omega)\times \dot {H}^{\beta-1}(\Omega)$ and $F\in L^{\tilde q'}((0,T); L^{\tilde r'}(\Omega))$.
 Then
\begin{align*}
\|u\|_{L^q((0,T);L^r(\Omega))}&+\|u\|_{L^\infty((0,T); \dot H^\beta(\Omega))}+\|\partial_t u\|_{L^\infty((0,T); \dot H^{\beta-1}(\Omega))}\\&\quad\leq C_{\text{Str}}\left(\|u_0\|_{\dot{H}^{\beta}(\Omega)}+\|u_1\|_{\dot{H}^{\beta-1}(\Omega)}+\|F\|_{L^{\tilde q'}((0,T); L^{\tilde r'}(\Omega))}\right),
\end{align*}
with 
\[
\frac{1}{q}\leq\frac{3}{4}\left(\frac{1}{2}-\frac{1}{r}\right) \,\,\,\text{and}\,\,\, \beta=3 \left(\frac{1}{2}-\frac{1}{r}\right)-\frac{1}{q}=3 \left(\frac{1}{2}-\frac{1}{\tilde r'}\right)-\frac{1}{\tilde q'}+2.
\]
\end{corollary}
\begin{proof}
The Duhamel formula yields
\[
u(t,.)=\cos(t\sqrt{-\Delta})u_0+\frac{\sin(t\sqrt{-\Delta})}{\sqrt{-\Delta}}u_1+\int_0^t\frac{\sin((t-s)\sqrt{-\Delta})}{\sqrt{-\Delta}}F(s,.)ds.
\]
The contribution of $(u_0, u_1)$ follows from Theorem \ref{STR}. It remains to prove the bounds on $(u,\partial_t u)$ in $\dot H^\beta(\Omega)\times \dot H^{\beta-1}(\Omega)$ and that
\[
\left\|\int_0^t\frac{\sin((t-s)\sqrt{-\Delta})}{\sqrt{-\Delta}}F(s,.)ds\right\|_{L^q((0,T);L^r(\Omega))}\lesssim \|F\|_{L^{\tilde q'}((0,T); L^{\tilde r'}(\Omega))}.
\]
By the  Christ-Kiselev lemma, it suffices to show that for $q>\tilde q',$
\begin{align}\label{CKL}
\left\|\int_0^T\frac{\sin((t-s)\sqrt{-\Delta})}{\sqrt{-\Delta}}F(s,.)ds\right\|_{L^q((0,T);L^r(\Omega))}\lesssim \|F\|_{L^{\tilde q'}((0,T); L^{\tilde r'}(\Omega))}.
\end{align}
Now, let $U(t)=e^{it\sqrt{-\Delta}}: L^2\to L^2$ be the half wave operator. Recall that
\[
\cos(t\sqrt{-\Delta})=\frac{U(t)+U(-t)}{2},\quad\sin(t\sqrt{-\Delta})=\frac{U(t)-U(-t)}{2i}.
\]
It follows from Theorem \ref{STR} that
\[
\|U(t)f\|_{L^q((0,T);L^r(\Omega))}\lesssim \|f\|_{\dot H^\beta(\Omega)}
\]
holds for all $(q,r,\beta)$ satisfying \eqref{ADM} and \eqref{ADM1}. For $\beta\in\mathbb{R}$ and $(q,r)$ satisfying satisfying \eqref{ADM} and \eqref{ADM1}, we define the operator $T_\beta$ by
\begin{align*}
T_\beta: L^2(\Omega)&\longrightarrow L^q((0,T);L^r(\Omega))\\
f&\longmapsto T_\beta(f)=(\sqrt{-\Delta})^{-\beta}e^{it\sqrt{-\Delta}}f.
\end{align*}
By the duality, it follows that the operator $T_{1-\beta}^*$ defined by
\begin{align*}
T_{1-\beta}^*: L^{\tilde q'}((0,T);L^{\tilde r'}(\Omega))&\longrightarrow L^2(\Omega)\\
F(s,.)&\longmapsto T_{1-\beta}^*(F(s,.))=\int_0^T(\sqrt{-\Delta})^{-(1-\beta)}e^{-is\sqrt{-\Delta}}F(s,.)ds,
\end{align*}
 where  $1-\beta=3\left(\frac{1}{2}-\frac{1}{\tilde r}\right)-\frac{1}{\tilde q}$. Hence, we get
 \begin{align*}
 \left\|\int_0^T U(t)U^*(s)(-\Delta)^{-\frac{1}{2}}F(s,.)ds\right\|_{L^q((0,T);L^r(\Omega))}&=\left\|T_\beta T_{1-\beta}^*F\right\|_{L^q((0,T);L^r(\Omega))}\\&\lesssim \|F\|_{ L^{\tilde q'}((0,T);L^{\tilde r'}(\Omega))},
 \end{align*}
 where $(q,r)$ and $(\tilde q,\tilde r)$ satisfy  \eqref{ADM} and \eqref{ADM1} as a consequence of $\beta=3 \left(\frac{1}{2}-\frac{1}{r}\right)-\frac{1}{q}$ and  $1-\beta=3\left(\frac{1}{2}-\frac{1}{\tilde r}\right)-\frac{1}{\tilde q}$. But
 \[
 \int_0^TU(t)U^*(s)(-\Delta)^{-\frac{1}{2}}F(s,.)ds=\int_0^T\frac{\sin((t-s)\sqrt{-\Delta})}{\sqrt{-\Delta}}F(s,.)ds,
 \]
 it follows that \eqref{CKL} holds. Observe that for all $(q,r)$ and $(\tilde q,\tilde r)$ satisfy  \eqref{ADM} and \eqref{ADM1}, we must have $q>\tilde q'$.
 
 To this end, we consider the bound on \[\|u(t,.)\|_{L^\infty((0,T); \dot H^\beta(\Omega))}+\|\partial_t u(t,.)\|_{L^\infty((0,T); \dot H^{\beta-1}(\Omega))}.\] We have
 \begin{align*}
 \|u(t,.)\|_{\dot H^\beta(\Omega)}+\|\partial_tu(t,.)\|_{\dot H^{\beta-1}(\Omega)}&=\left\|\cos(t\sqrt{-\Delta})u_0+\frac{\sin(t\sqrt{-\Delta})}{\sqrt{-\Delta}}u_1\right\|_{\dot H^\beta(\Omega)}\\&+\left\|(-\sqrt{-\Delta})\sin(t\sqrt{-\Delta})u_0+\cos(t\sqrt{-\Delta})u_1\right\|_{\dot H^{\beta-1}(\Omega)}\\
 &\leq 2(\|u_0\|_{\dot H^\beta(\Omega)}+\|u_1\|_{\dot H^{\beta-1}(\Omega)}).
 \end{align*}
 It  reduces to showing that it is uniformly in $t$ such that
 \begin{align*}
 \left\|\int_0^t\frac{\sin((t-s)\sqrt{-\Delta})}{\sqrt{-\Delta}}F(s,.)ds\right\|_{\dot H^\beta(\Omega)}+&\left\|\int_0^t\cos((t-s)\sqrt{-\Delta})F(u(s,.))ds\right\|_{\dot H^{\beta-1}(\Omega)}\\&\lesssim \|F\|_{L^{\tilde q'}((0,T); L^{\tilde r'}(\Omega))}.
 \end{align*}
 or equivalently,
 \begin{equation}\label{INH}
  \begin{aligned}
 \left\|\int_0^t\frac{\sin((t-s)\sqrt{-\Delta})}{(\sqrt{-\Delta})^{1-\beta}}F(s,.)ds\right\|_{L^2(\Omega)}+&\left\|\int_0^t\frac{\cos((t-s)\sqrt{-\Delta})}{(\sqrt{-\Delta})^{1-\beta}}F(s,.)ds\right\|_{L^2(\Omega)}\\&\lesssim \|F\|_{L^{\tilde q'}((0,T); L^{\tilde r'}(\Omega))}.
 \end{aligned}
 \end{equation}
 But as consequence of the boundedness of the adjoint operator $T_{1-\beta}^*$, we obtain
 \[
 \left\|\int_0^T\frac{e^{i(t-s)\sqrt{-\Delta}}}{(\sqrt{-\Delta})^{1-\beta}}F(s,.)ds\right\|_{L^2(\Omega)}\lesssim \|F\|_{L^{\tilde q'}((0,T); L^{\tilde r'}(\Omega))}.
 \]
 and the Christ-Kiselev lemma yields
  \[
 \left\|\int_0^t\frac{e^{i(t-s)\sqrt{-\Delta}}}{(\sqrt{-\Delta})^{1-\beta}}F(s,.)ds\right\|_{L^2(\Omega)}\lesssim \|F\|_{L^{\tilde q'}((0,T); L^{\tilde r'}(\Omega))}.
 \]
 and \eqref{INH} follows from Euler formula. This completes the proof.
\end{proof}

\section{Well-posedness}
\subsection{Local Existence}
In this section, we establish the local existence for the solution to the wave equation \eqref{NLW} by applying the Strichartz estimates in Theorem \ref{STR} and the fixed point argument.

 Theorem \ref{LOC} follows immediately from the following result.
\begin{theorem}
Let $(\Omega,\Delta)$ be defined as before. For any $(u_0,u_1)\in H_0^1(\Omega)\times L^2(\Omega)$. Then there exists $T=T(\|(u_0,u_1)\|_{H_0^1(\Omega)\times L^2(\Omega)})>0$ such that the energy-critical wave equation \eqref{NLW} is local wellposed in $(0,T)$ and the unique solution $u$ satisfies
\[
u\in C^0((0,T); H_0^1(\Omega))\cap C^1((0,T); L^2(\Omega))\cap L^5((0,T); L^{10}(\Omega)).
\]
Moreover, if $\|(u_0,u_1)\|_{H_0^1(\Omega)\times L^2(\Omega)}<\delta$ for small enough $\delta$, there exists a global solution $u$.
\end{theorem}

\begin{proof}
We apply the Banach fixed point argument to prove this result. Let a space $X_T$ be defined by
\[
X_T=C^0((0,T); H_0^1(\Omega))\cap C^1((0,T); L^2(\Omega))\cap L^5((0,T); L^{10}(\Omega)).
\]
We define a norm $\|.\|_T$ by
\[
\|u\|_{X_T}=\sup_{0<t<T}\|(u,\partial_t u)\|_{H_0^1(\Omega)\times L^2(\Omega)}+\|u\|_{L^5((0,T); L^{10}(\Omega))}.
\]
We remark that the space $(X_T, \|.\|_T)$ is a Banach space.  For any small constant $\varepsilon>0$. We introduce a fixed-point space
\begin{align*}
B_T=\{u\in X_T: \|(u,\partial_tu)\|_{H_0^1(\Omega)\times L^2(\Omega)}\leq 2C_{\text{Str}}\delta, \|u\|_{L^5((0,T); L^{10}(\Omega))}\leq 2C_{\text{Str}}\varepsilon\}
\end{align*}
with the metric
\[d(u,v)=\sup_{0<t<T}\|(u-v,\partial_tu-\partial_t v)\|_{H_0^1(\Omega)\times L^2(\Omega)}+\|u-v\|_{L^5((0,T); L^{10}(\Omega))}.\]
 Consider the solution map
\begin{align}\label{FPM}
\Phi(u)&=\cos(t\sqrt{-\Delta})u_0+\frac{\sin(t\sqrt{-\Delta})}{\sqrt{-\Delta}}u_1+\int_0^t\frac{\sin((t-s)\sqrt{-\Delta})}{\sqrt{-\Delta}}u^5(s,.)ds.
\end{align}
 Let denote 
\[
u_{\text{hom}}:=\cos(t\sqrt{-\Delta})u_0+\frac{\sin(t\sqrt{-\Delta})}{\sqrt{-\Delta}}u_1
\]
and 
\[
u_{\text{inh}}:=\int_0^t\frac{\sin((t-s)\sqrt{-\Delta})}{\sqrt{-\Delta}}u^5(s,.)ds.
\]

It follows from the Strichartz estimates in Theorem \ref{STR} with the admissible triplet $q=5, r=10,$ and $\beta=1$ that
 \begin{align*}
 \|u\|_{L^5((0,T); L^{10}(\Omega))}+\|u\|_{L^\infty((0,T); H_0^1(\Omega))}&+\|\partial_tu\|_{L^\infty((0,T); L^2(\Omega))}\\&\leq C_{\text{Str}} (\|u_0\|_{H_0^1(\Omega)}+\|u_1\|_{L^2(\Omega)}).
 \end{align*}
 This yields
 \[
 \left\|u_{\text{hom}}\right\|_{L^5((0,T); L^{10}(\Omega))}\leq C_{\text{Str}} (\|u_0\|_{H_0^1(\Omega)}+\|u_1\|_{L^2(\Omega)}).
 \]
Therefore, if the  norm of initial data $\|(u_0,u_1)\|_{H_0^1(\Omega)\times L^2(\Omega)}<\delta$ for small enough $\delta$, then we have
 \[
 \left\|u_{\text{hom}}\right\|_{L^5((0,T); L^{10}(\Omega))}\leq C_{\text{Str}}\varepsilon
 \]
 holds for $T=\infty$; otherwise, the inequality holds for some small $T>0$ as a consequence of the dominated convergence theorem. It follows from Corollary \ref{INS} with admissible triplet $q=5, r=10,$ and $\beta=1$ (so that $\tilde q'=1$ and $\tilde r'=2$) that
 \[
 \|u_{\text{inh}}\|_{L^5((0,T); L^{10}(\Omega))}\leq C_{\text{Str}}\|u^5\|_{L^1((0,T); L^2(\Omega))}.
 \] 
 We need to show that the operator $\Phi$ is well-defined on $B_T$ and is a contraction map under the metric $d$. To do this, let $u\in B_T$ with $0<\varepsilon\ll 1$. Then, by Strichartz estimates we have
 \begin{align*}
 \|\Phi(u)\|_{L^5((0,T); L^{10}(\Omega))}&\leq\|u_{\text{hom}}\|_{L^5((0,T); L^{10}(\Omega))}+ \|u_{\text{inh}}\|_{L^5((0,T); L^{10}(\Omega))}\\
& \leq C_{\text{Str}}\varepsilon+C_{\text{Str}}\|u^5\|_{L^1((0,T); L^2(\Omega))}\\
& \leq2C_{\text{Str}}\varepsilon,
 \end{align*}
  and
 \begin{align*}
 \sup_{0<t<T}\|(\Phi(u),\partial_t\Phi(u))\|_{H_0^1(\Omega)\times L^2(\Omega)}&\leq C_{\text{Str}}\|(u_0,u_1)\|_{H_0^1(\Omega)\times L^2(\Omega)}+C_{\text{Str}}\|u^5\|_{L^1((0,T); L^2(\Omega))}\\
 &\leq 2C_{\text{Str}}\|(u_0,u_1)\|_{H_0^1(\Omega)\times L^2(\Omega)}\\
 &\leq 2C_{\text{Str}}\delta,
 \end{align*}
 for small $T>0$. Therefore, $\Phi(u)\in B_T$. Now, we prove that $\Phi$ is a contraction map. Let $v_1, v_2\in B_T$. It follows from the Strichartz estimates and by choosing $\varepsilon$ sufficiently small, we get
 \begin{align*}
 \|\Phi(v_1)-\Phi(v_2)\|_{L^5((0,T); L^{10}(\Omega))}&\leq C_{\text{Str}}\|v_1^5-v_2^5\|_{L^1((0,T); L^2(\Omega))}\\
 &\leq C_{\text{Str}}\|v_1-v_2\|_{L^5((0,T); L^{10}(\Omega))}\\&\quad\times\big(\|v_1\|_{L^5((0,T); L^{10}(\Omega))}^4+\|v_2\|_{L^5((0,T); L^{10}(\Omega))}^4\big)\\
 &\leq\tilde C\varepsilon^4 \|v_1-v_2\|_{L^5((0,T); L^{10}(\Omega))},
 \end{align*}
 and 
 \begin{align*}
  \sup_{0<t<T}\|(\Phi(v_1)-\Phi(v_2),\partial_t\Phi(v_1)-\partial_t\Phi(v_2))&\|_{H_0^1(\Omega)\times L^2(\Omega)}\\
  &\leq C_{\text{Str}}\|v_1^5-v_2^5\|_{L^1((0,T); L^2(\Omega))}\\
   &\leq\tilde C\varepsilon^4 \|v_1-v_2\|_{L^5((0,T); L^{10}(\Omega))}.
 \end{align*}
 We combine these two estimates to obtain
 \begin{align*}
 d(\Phi(v_1)-\Phi(v_2))\leq \tilde C\varepsilon^4 \|v_1-v_2\|_{L^5((0,T); L^{10}(\Omega))}\leq\frac{1}{2}d(v_1,v_2).
 \end{align*}
 We conclude from the fixed point theorem that there is a unique solution $u$ to \eqref{NLW} on $(0, T)\times \Omega$. Hence, if $\delta$ is small enough, we get the global solution; otherwise, we have a local solution.
\end{proof}
\subsection{Global Existence}
In this section, we follow the ideas in \cite{BLP} as well as  \cite{ShS,ShS1} to obtain the global solution to \eqref{NLW} for arbitrary finite energy data in the context of our cylindrical convex domains as stated in Theorem \ref{GWP}. The key ingredients are the following stronger version of Strichartz estimates, trace estimates and the nonconcentration of nonlinear effect in a small light cones.
\begin{theorem}\label{SSTR}
Let $(\Omega, \Delta)$ be defined as before.
If $u, F$ satisfy
\begin{align*}
&(\partial_t^2-\Delta)u=F,\,\,\,u_{|t=0}=u_0,\,\,\, \partial_t u_{|t=0}=u_1,\,\,\,u_{|x=0}=0,
\end{align*}
then 
\begin{align*}
\|u\|_{L^5\big((0,1);W_0^{\frac{3}{10},5}(\Omega)\big)}&+\|u\|_{C^0((0,1);  H_0^1(\Omega))}+\|\partial_t u\|_{C^0((0,1); L^2(\Omega))}\\&\leq C_{\text{Str}}\left(\|u_0\|_{H_0^{1}(\Omega)}+\|u_1\|_{L^2(\Omega)}+\|F\|_{L^{\frac{5}{4}}\big((0,1); W^{\frac{7}{10},\frac{5}{4}}(\Omega)\big)}\right).
\end{align*}
\end{theorem}
\begin{proof}
We follow the streamline of the proof of Proposition 3.1 in \cite{BLP} with some modifications for the homogeneous part. We have by Duhamel formula
\[
u(t,.)=\cos(t\sqrt{-\Delta})u_0+\frac{\sin(t\sqrt{-\Delta})}{\sqrt{-\Delta}}u_1+\int_0^t\frac{\sin((t-s)\sqrt{-\Delta})}{\sqrt{-\Delta}}F(s,.)ds.
\]
 The contribution of $(u_0,u_1)$ follows from the Strichartz estimates in Theorem \ref{STR} with the admissible triplet $q=5, r=5,$ and $\beta=\frac{7}{10}$ that
 \begin{equation}\label{STR5}
 \left\|e^{it\sqrt{-\Delta}}f\right\|_{L^5((0,1)\times\Omega)}\lesssim \|f\|_{\dot H^{\frac{7}{10}}(\Omega)}.
 \end{equation}
 Then if we apply this inequality to $\Delta f$ and we use the $L^p$ elliptic regularity, we get
   \begin{equation}\label{LPR}
 \left\|e^{it\sqrt{-\Delta}}f\right\|_{L^5((0,1); W^{2,5}(\Omega)\cap W_0^{1,5}(\Omega))}\lesssim \|f\|_{\dot H^{\frac{27}{10}}(\Omega)}.
 \end{equation}
 Consequently, the interpolation between \eqref{STR5} and \eqref{LPR} gives
 \[
  \left\|e^{it\sqrt{-\Delta}}f\right\|_{L^5\big((0,1); W_0^{\frac{3}{10},5}(\Omega)\big)}\lesssim \|f\|_{\dot H^{1}(\Omega)}.
  \]
  We can conclude that
  \[
  \left\|\cos(t\sqrt{-\Delta})u_0+\frac{\sin(t\sqrt{-\Delta})}{\sqrt{-\Delta}}u_1\right\|_{L^5\big((0,1); W_0^{\frac{3}{10},5}(\Omega)\big)}\lesssim \|u_0\|_{H_0^1(\Omega)}+\|u_1\|_{L^2(\Omega)}.
  \]
  The bound for \[\|u\|_{C^0((0,1);  H_0^1(\Omega))}+\|\partial_t u\|_{C^0((0,1); L^2(\Omega))}\] is similar to that of in Theorem \ref{STR}. Finally, the bound of 
  \[
  \int_0^t \frac{e^{i(t-s)\sqrt{-\Delta}}}{\sqrt{-\Delta}}F(s,.)ds
  \]
   follows that same line as in the proof of Proposition 3.1 in \cite{BLP} by $TT^*$ argument and Christ-Kiselev lemma.
\end{proof}
Notice that one has the control for the nonlinear term (see \cite{BLP}) by the following estimate
\[
\|u^5\|_{L^{\frac{5}{4}}\big((0,1); W^{\frac{7}{10},\frac{5}{4}}(\Omega)\big)}\lesssim \|u\|_{L^5((0,1); L^{10}(\Omega))}^4\|u\|_{L^\infty((0,1); L^6(\Omega))}^{\frac{3}{10}}\|u\|_{L^\infty((0,1); H^1(\Omega))}^{\frac{7}{10}}.
\]
The idea now is to localize these estimates on small light cones and use the fact that the $L_t^\infty(L_x^6)$ norm is small in such a small cones (see \cite{BLP,ShS,ShS1} for details). Let us point out the key results to obtain the global existence to the \eqref{NLW} as follows.
\begin{proposition}[$L^6$-nonconcentration]
Let $x_0\in\overline{\Omega}$. Then for any solution $u$ to \eqref{NLW} in the space $X_{<t_0}=C^0((0, t_0); H_0^1(\Omega))\cap L_{loc}^5((0,t_0); L^{10}(\Omega))\times C^0((0,t_0); L^2(\Omega))$, there holds
\[
\lim_{t\to t_0^-}\int_{x\in\Omega\cap\{|x-x_0|<t-t_0\}}u^6(t,x)dx=0
\]
\end{proposition}
For the proof of this result see \cite{BLP}. We point out here three important ingredients. 
\begin{itemize}
\item The normal derivative estimate 
\[
\|\partial_\nu u\|_{L^2((0,t_0)\times \Omega)}\leq C E(u)^{\frac{1}{2}}
\]
is satisfied uniformly for $0<T<t_0$, where $\partial_\nu u$ is the trace to the boundary of the outward unit normal of $u$. This follows from the estimate (see Proposition 3.2 in \cite{BLP})
\[
\int_0^T\int_\Omega\big[(\partial_t^2-\Delta)Zu(t,x)-Z(\partial_t^2-\Delta)u(t,x)\big]u(t,x)dxdt\lesssim E(u)
\]
uniformly in $T$, where $Z$ is a smooth vector filed on $\Omega$ which coincides with $\partial_\nu$ on $\partial\Omega$. Notice that such an estimate follows by using integration by parts and the energy conservation.
\item The flux identity derived in \cite{BLP}
\begin{align*}
\int_{x\in\Omega,|x|<-T}&\left(\frac{1}{2}|\partial_t u|^2+\frac{1}{2}|\nabla u|^2+\frac{1}{6}|u|^6\right)(x,T)dx+\text{Flux}(u, M_S^T)\\
&=\int_{x\in\Omega,|x|<-S}\left(\frac{1}{2}|\partial_t u|^2+\frac{1}{2}|\nabla u|^2+\frac{1}{6}|u|^6\right)(x,S)dx=E_{loc}(S).
\end{align*}
Here,
\[
\text{Flux}(u, M_S^T):=\int_{M_S^T}\langle e(u),\nu\rangle d\sigma(t,x)
\]
where $M_S^T:=\{x: |x|=-t, S<t<T\}$, $\nu$ the unit outward normal to $M_S^T$, $d\sigma(t,x)$ the induced measure on $M_S^T$, and the vector field $e$ is given by
\[
e(u)=\left(\frac{1}{2}|\partial_t u|^2+\frac{1}{2}|\nabla u|^2+\frac{1}{6}|u|^6, -\partial_tu\nabla u\right).
\]
The authors in \cite{BLP} showed that
\[
\lim_{S\to 0^-}\text{Flux}(u, M_S^0)=0.
\]
\item A Morawetz type inequality is formally derived in \cite{BLP} and integrate the identity over $K_S^T=\{(t,x); |x|<-t, S<t<T\}\cap\Omega$ combined with H\"{o}lder's inequality and the conservation of energy. More precisely, they proved that
\[
\int_{x\in \Omega,|x|<-S}\frac{|u|^6}{6}(S,x)dx\leq |S|E(u)+C\text{Flux}(u, M_S^0)+C\text{Flux}(u, M_S^0)^{\frac{1}{3}}\underset{S\to 0^-}{\longrightarrow 0}.
\]
\end{itemize}
Finally, the next proposition shows a localizing space-time estimates.
\begin{proposition}[see \cite{BLP}]
For any $\varepsilon>0$, there exists $t<0$ such that
\[
\|u\|_{(L^5; L^{10})(K_t^0)}<\varepsilon,
\]
where
\[
\|u\|_{(L^5; L^{10})(K_t^0)}=\left(\int_{s=t}^0\left(\int_{\{|x|<-s\}\cap\Omega}|u|^{10}(s,x)dx\right)^{\frac{5}{10}}ds\right)^{\frac{1}{5}}.
\]
\end{proposition}
With these results in hand, to prove the global existence it is sufficient to show that the local solution $u$ satisfies $u\in L^5((s,0); L^{10}(\Omega))$ since this allows us to show that $\lim_{s'\to 0^-}(u,\partial_s u)(s',.)$ exists in $H_0^1(\Omega)\times L^2(\Omega)$ via Duhamel formula and the Strichartz estimates in Theorem \ref{SSTR} together with conservation of energy.
\subsection*{Acknowledgement}
The author would like to thank the referees for their insightful remarks and comments on this paper. The author is grateful for the hospitality of Link\"{o}ping University where the project was initiated under a grant of the Abel Visiting Scholar Program of International Mathematical Union.
\bibliographystyle{amsplain}
\bibliography{mybib1}
\end{document}